\documentclass[reqno]{amsart}
\usepackage{amsmath,amssymb}
\usepackage{hyperref}

\usepackage[T1]{fontenc}
\usepackage[utf8]{inputenc}
\usepackage[english]{babel}
\usepackage{csquotes}

\newcommand*{\mailto}[1]{\href{mailto:#1}{\nolinkurl{#1}}}

\usepackage{amsfonts}
\usepackage{amssymb}
\usepackage{amsthm}
\usepackage{enumitem}

\newcommand*{\oo}{\infty}
\newcommand*{\N}{\mathbb{N}}
\newcommand*{\Z}{\mathbb{Z}}

\newcommand*{\R}{\mathbb{R}}

\newcommand*{\eps}{\varepsilon}

\newcommand*{\sumz}{\sum_{n \in \Z}}

\newtheorem{thm}{Theorem}[section]
\newtheorem{lem}[thm]{Lemma}
\newtheorem{rem}[thm]{Remark}
\numberwithin{equation}{section}

\begin{document}

\title[Global existence versus blow-up for FPUT lattices]{Global existence versus blow-up for Fermi--Pasta--Ulam--Tsingou lattices}

\author[G. Teschl]{Gerald Teschl}
\address{Faculty of Mathematics\\ University of Vienna\\
Oskar-Morgenstern-Platz 1\\ 1090 Wien\\ Austria}
\email{\mailto{Gerald.Teschl@univie.ac.at}}
\urladdr{\url{https://www.mat.univie.ac.at/~gerald/}}

\author[A. Zand]{Anton Zand}
\address{Faculty of Mathematics\\ University of Vienna\\
Oskar-Morgenstern-Platz 1\\ 1090 Wien\\ Austria}
\email{\mailto{anton.zand@univie.ac.at}}


\begin{abstract} 
This paper investigates the global behavior of infinite nonlinear lattices. In particular, for the Fermi--Pasta--Ulam--Tsingou-$\alpha$ lattice (FPUT-$\alpha$), we show global existence for small initial data via a region-based approach, and prove blow-up for negative energy and some positive energy, given appropriate conditions on the initial momentum and dispersion.

Moreover, we extend these results to the nonlinear Klein-Gordon (DNKG) equation and nonlinear wave (DNW) equation.
\end{abstract}

\maketitle

\section{Introduction}

The study of nonlinear lattices has both a rich history, including the pioneering work in numerics by Fermi, Pasta, Ulam and Tsingou \cite{fpu}, as well as a wide range of applications in mechanics, optics or acoustics. Despite the vast amount of literature on FPUT lattices, there are still questions about the long-time behavior of certain models, such as the ones posed by Pankov in \cite{pank}. In particular, the aim of the present paper is to provide a global existence result for the FPUT-$\alpha$ model with small positive energy, providing a solution of
Problem~1.1 from Pankov \cite{pank}. Moreover, we complement this result by establishing blow-up results for initial conditions with negative energy and certain initial conditions
with positive energy. To this end we adapt a Levine--Glassey-type concavity argument from PDE-theory. We also provide examples of initial conditions satisfying our assumptions.

Finally, we extend these results to the nonlinear Klein-Gordon (DNKG) equation for the focusing case with even exponent and its massless limit, the nonlinear wave (DNW) equation.
Further results can be found in \cite{kar, zand}.

\section{The Fermi--Pasta--Ulam--Tsingou Lattice}

We consider an infinite one-dimensional lattice assuming nearest neighbor interaction and unitary masses. The differential equation for the particle displacement $q_n(t)$ of such a system reads 
\begin{equation}\label{C1}
    \ddot q_n(t) = V'(q_{n+1}-q_{n})-V'(q_{n}-q_{n-1})
\end{equation}
with some potential $V\in C^2(\R)$. We assume $V(0)=V'(0)=0$ as usual. We also introduce the momentum $p=\dot q$, the differences of displacements $r_n=q_{n+1}-q_n$ and abbreviate $\mathfrak{D}(t) := \|r(t)\|_\oo$. Concerning the phase space, we take initial data in $\ell^2(\Z)\times\ell^2(\Z)$, which ensures $\lim_{n \to \oo} q_n = 0$. Existence and uniqueness of local solutions is standard. It is also well known that the total energy of such a system, given by the Hamiltonian
\begin{equation}\label{H1}
    \mathcal{H}(p,q):=E_k + E_p=\sum_{n \in \Z} \frac{p_n^2}{2} +\sum_{n \in \Z} V(q_{n+1}-q_n),
\end{equation}
is a conserved quantity (we will hence write $\mathcal{H}(p,q)(t)=\mathcal{H}$ for brevity). If $V \ge 0$ and $\mathcal{H}<\oo$ for all $t$, the bound
\begin{equation}\label{SOL}
    \|p(t)\|_2^2 \leq 2\mathcal{H}, \quad \|q(t)\|_2\leq \|q(0)\|_2 +\sqrt{2\mathcal{H}}t
\end{equation}
ensures that solutions are global by standard theory. We refer to \cite{pank} or \cite[Ch.~9.6]{fa} for details.

The FPUT-$\alpha$ model arises from setting $V(r):=\frac{1}{2}r^2+\frac{\alpha}{3}r^3$. In contradistinction to the FPUT-$\beta$ lattice or the Toda lattice (see \cite{jaco,toda}),
$V$ does not satisfy the nonnegativity assumption and hence the above argument does not apply. As our first result we establish that we have global existence for small initial data:

\begin{thm}[Global Existence for Small Initial Data]\label{GLO}
    Let $\alpha \ne 0$, $0<A<\frac{3}{2}$ and define $B$ by $\frac{1}{B}=\frac{1}{2}-\frac{A}{3}$. For initial data in $\ell^2\times\ell^2$, solutions of the FPUT-$\alpha$ equation satisfying $\mathfrak{D}(0) < \frac{A}{|\alpha|} \quad \text{and} \quad 0 < \mathcal{H} < \frac{A^2}{B \alpha^2}$ are global and remain within the region $\mathfrak{D}(t) < \frac{A}{|\alpha|}$.
\end{thm}

\begin{proof}
    $A$ and $B$ are chosen such that $0 \leq \frac{r^2}{B} \leq V(r) \leq r^2$ holds for as long as $\mathfrak{D}(t) \leq \frac{A}{|\alpha|}$. Thus, under this restriction, $\mathcal{H}$ is finite, hence well defined and conserved, by $\infty > \|p\|_2^2 + 4\|q\|_2^2 \geq \sum_{n \in \Z} \frac{p_n^2}{2} + V(q_{n+1}-q_n) = \mathcal{H}$ (using {$V(r) \le r^2$} and {$(a-b)^2 \le 2a^2+2b^2$}). On the other hand, $\nonumber \frac{A^2}{B\alpha^2} > \mathcal{H} \geq \frac{\mathfrak{D}(t)^2}{B}$ forces the solution to remain within the region. Finally, \eqref{SOL} finishes the proof since $V_\alpha \ge 0$ is guaranteed.
\end{proof}

To show blow-up, we adapt a concavity argument first used in the PDE context by Levine in \cite{levi}. The specific structure used here follows Glassey \cite{glass}. We first define the functional $M(t) := \frac{1}{2}\sum_{n\in\mathbb{Z}} q_n(t)^2$ and prove a virial-type inequality.

\begin{lem}[Virial Inequality]\label{VIR}
Consider the FPUT-$\alpha$ model. We have 
\begin{equation}
(1+\gamma)\dot{M}^2 - M\ddot{M} \le 3\mathcal{H}M, \quad \gamma := \frac{1}{4}.
\end{equation}
Further, in the case $\mathcal{H}<0$, there exists $\eps$ such that $M(t)\ge \eps >0$ for all $t$.
\end{lem}

\begin{proof}
Differentiating $M(t)$ twice and using summation by parts along with equation \eqref{C1} gives
\begin{equation}
\ddot{M} = 2E_k - \sum_{n\in\mathbb{Z}} r_n V'(r_n).
\end{equation}
Observe $r V'(r) = 3V(r) - \frac{1}{2}r^2$ and $\sumz V(r_n) = \mathcal{H} - E_k$, hence we may substitute and obtain
\begin{equation}
\ddot{M} = 5E_k - 3\mathcal{H} + \frac{1}{2} \sumz r_n^2 \ge 5E_k - 3\mathcal{H}.
\end{equation}
The Cauchy--Schwarz inequality gives $\dot{M}^2 = (\sum_n q_n \dot{q}_n)^2 \le \|q\|_2^2\|p\|_2^2 = 4M E_k$, hence $E_k \ge \frac{\dot{M}^2}{4M}$, so finally
\begin{equation}
\ddot{M} \ge \frac{5\dot{M}^2}{4M} - 3\mathcal{H} = (1+\gamma)\frac{\dot{M}^2}{M} - 3\mathcal{H},
\end{equation}
and we conclude $(1+\gamma)\dot{M}^2 - M\ddot{M} \le 3\mathcal{H}M$.

For the second claim, observe $E_k \ge 0$, so $E_p \le \mathcal{H} < 0$. Achieving $E_p < 0$ requires at least one $n$ such that $V(r_n) < 0$, implying $|r_n| > \frac{3}{2|\alpha|}$ for some $n$. Minimizing $q_{n+1}^2 + q_n^2$ subject to $|q_{n+1}| + |q_n| \ge |r_n| > \frac{3}{2|\alpha|}$ shows
\begin{equation}
M(t) \ge \frac{q_{n+1}(t)^2 + q_n(t)^2}{2} \ge \frac{9}{16\alpha^2} > 0.
\end{equation}
\end{proof}

\begin{rem}
For $M>0$, the Cauchy--Schwarz defect $E_k-\frac{\dot{M}^2}{4M}$ is nonnegative, and retaining it gives indeed the sharper identity
\[
    (1+\gamma)\dot{M}^2-M\ddot{M}
    =3\mathcal{H}M-\frac{M}{2}\sum_n r_n^2
      -5M\left(E_k-\frac{\dot{M}^2}{4M}\right).
\]
\end{rem}

\begin{thm}[Blow-up for Negative Energy]\label{E-}
Solutions of the FPUT-$\alpha$ model with $\mathcal{H} < 0$ are not global in time.
\end{thm}

\begin{proof}
Define $N(t) := M(t)^{-\gamma}$. Differentiating twice yields
\begin{equation}\label{.VIR}
\ddot{N} = \gamma M^{-\gamma-2} \left( (\gamma+1)\dot{M}^2 - M\ddot{M} \right).
\end{equation}
By Lemma \ref{VIR}, $M\ge\eps>0$ and the bracket in \eqref{.VIR} is strictly negative, so $N$ is strictly concave. Also, the estimate $\ddot M\ge -3\mathcal{H}>0$ implies that $\dot M$ eventually becomes positive if the solution persists. At such a time $t_1$, $\dot N(t_1)<0$. Strict concavity then forces $N$ to reach zero within finite forward time; equivalently, $M$ becomes unbounded and the solution cannot remain in $\ell^2$.
\end{proof}

\begin{thm}[Blow-Up for Positive Energy]\label{E+}
Solutions of the FPUT-$\alpha$ model with $\mathcal{H} > 0$ and initial momentum subject to $i)$ $\dot{M}(0) > 0$ and $ii)$ $\dot{M}(0)^2 > 4\mathcal{H}M(0)$ are not global in forward time. Moreover, their maximal forward existence time satisfies
\begin{equation*}
    T_+ \leq \frac{M(0)^{-\gamma}}{\sqrt{\big(-\gamma M(0)^{-\gamma-1}\dot M(0)\big)^2 - \gamma \mathcal{H}M(0)^{-6\gamma}}}, \qquad \gamma=\frac14.
\end{equation*}
\end{thm}

\begin{proof}
As before we define $N(t) := M(t)^{-\gamma}$. 
By Lemma \ref{VIR} and \eqref{.VIR},
\begin{equation}
\ddot{N} = \gamma M^{-\gamma-2} \left( (\gamma+1)\dot{M}^2 - M\ddot{M} \right) \le \gamma M^{-\gamma -2} \big(3M\mathcal{H}\big) = 3\gamma N^5 \mathcal{H}.
\end{equation}
Observe that $\dot N(t) < 0$ at $0$ by $i)$, and by continuity on some interval $I := [0,T)$, which we take to be maximal. Thus multiplying by $2\dot N$ and integrating from $0$ to $t$ yields
\begin{equation}
\dot{N}(t)^2 \ge \gamma\mathcal{H}N(t)^6 - \gamma\mathcal{H}N(0)^6 + \dot{N}(0)^2 \ge \zeta > 0,
\end{equation}
where $\zeta := \dot{N}(0)^2 - \gamma\mathcal{H}N(0)^6 > 0$ by $ii)$. Consequently, $\dot{N}(t) \le -\sqrt{\zeta} < 0$ for all $t \in I$, and the lower bound on $\dot N(t)^2$ ensures that $I$ includes the whole forward existence interval. Integrating once more gives $0\leq N(t)\leq N(0)-\sqrt{\zeta}\,t$, so $N$ reaches zero within time $N(0)/\sqrt{\zeta}$. Thus $M(t)\to\infty$ at the maximal forward existence time, and
\[
    T_+\leq \frac{N(0)}{\sqrt{\zeta}}
    =\frac{M(0)^{-\gamma}}{\sqrt{\big(-\gamma M(0)^{-\gamma-1}\dot M(0)\big)^2-\gamma\mathcal{H}M(0)^{-6\gamma}}}.
\]
\end{proof}

\section{The Nonlinear-Klein--Gordon Lattice}

In this section we extend our results to the discrete Nonlinear-Klein--Gordon (DNKG) equation. Here, we consider a linear nearest-neighbor interaction with an additional nonlinear potential.
Such a system is given by
\begin{equation}\label{C2}
    \ddot q_n = (Aq)_n-V'(q_n),
\end{equation}
where specifically for the DNKG equation we have
\begin{equation}
    \begin{aligned}
        A&=a\Delta-m^2, \qquad (\Delta q)_n=q_{n+1}-2q_n+q_{n-1},\\
        V(r)&=-\frac{\lambda}{k+1}r^{k+1}, \qquad m,\lambda\in\R,\quad a\in\R_+,\quad k\in\N.
    \end{aligned}
\end{equation}
The limiting case $m=0$ is known as the discrete nonlinear wave (DNW) equation.

Again, the Hamiltonian
\begin{equation}\label{eq:ham2}
    \mathcal{H}(p,q):= E_k + E_l+E_p =\sum_{n \in \Z} \bigg(\frac{p_n^2}{2} - \frac{(A q)_nq_n}{2} + V(q_n)\bigg),
\end{equation}
is a conserved quantity for data in $\ell^2\times\ell^2$. For $k$ odd and $\lambda <0$, solutions are global by $V\ge0$ as already mentioned.
Hence we focus on the case where $k$ is even. For the special (cubic) case $k=2$ of the DNKG, global existence for small initial data is established by Pankov \cite[Thm.~1.9]{pank}.

The proof of blow-up is similar to the FPUT-$\alpha$ model. We define again $M(t) := \frac{1}{2}\sum_{n \in \Z}q_n(t)^2.$

\begin{lem}[Virial Inequality]\label{thm:lckg1}
    Consider the DNKG or DNW equation with $\lambda > 0$, $k \geq 2$ even. Then it holds that
    \begin{equation}
        (\gamma +1) \dot M^2 - M\ddot M \leq (k+1)\mathcal{H}M-(k-1)M E_l, \quad \gamma := \frac{k-1}{4}.
    \end{equation}
    If we additionally assume $\mathcal{H} < 0$, then there exists a constant $\eps >0$ such that $M(t) \geq \eps > 0$ for all $t$.
\end{lem}

\begin{proof}
    It is readily verified that $E_l\geq0$, thus $E_p\leq \mathcal{H}<0$, meaning $\frac{-\lambda}{k+1} \|q\|_{k+1}^{k+1}\leq \mathcal{H}$. We deduce $-\frac{k+1}{\lambda} \mathcal{H} \leq \|q\|_2^{k+1} = \big(2M\big)^{\frac{k+1}{2}}$ and further $M \geq \frac{1}{2} \big(-\frac{k+1}{\lambda} \mathcal{H}\big)^{\frac{2}{k+1}} > 0$.
    
    Differentiating $M$ twice and using the energy identity gives
    \begin{equation}\label{eq:dnkg-virial}
        \ddot M=(k+3)E_k+(k-1)E_l-(k+1)\mathcal{H}.
    \end{equation}
    Since $\dot M^2\leq 4ME_k$ and $4(\gamma+1)=k+3$, it follows that
    \[
        (\gamma+1)\dot M^2-M\ddot M
        \leq (k+1)\mathcal{H}M-(k-1)ME_l.\qedhere
    \]
\end{proof}

\begin{thm}[Blow-Up for Negative Energy]\label{DNKG}
    Solutions of the DNKG or the DNW equation with $\lambda > 0$, $k \geq 2$ even and $\mathcal{H} < 0$ are not global in time.
\end{thm}

\begin{proof}
    Set $N(t):=M(t)^{-\gamma}$ with $\gamma=(k-1)/4$. The virial inequality gives
    \[
        \ddot N \leq \gamma(k+1)\mathcal{H}M^{-\gamma-1}<0.
    \]
    Moreover, $\ddot M\geq -(k+1)\mathcal{H}>0$, so $\dot M$ eventually becomes positive if the solution persists. Then $\dot N<0$, and strict concavity forces $N$ to reach zero in finite forward time. Thus $M$ becomes unbounded and the solution cannot remain in $\ell^2$.
\end{proof}

\noindent As before, we also have a result for positive energy.

\begin{thm}[Blow-Up for Positive Energy]\label{DNKG-pe}
    Solutions of the DNKG or the DNW equation with $\lambda > 0$, $k \geq 2$ even and $\mathcal{H} > 0$ subject to $i)$ $\dot{M}(0) > 0$ and $ii)$ $\dot{M}(0)^2 > 4\mathcal{H}M(0)$ are not global in forward time. Moreover, with $\gamma=(k-1)/4$, their maximal forward existence time satisfies
    \begin{equation*}
        T_+ \leq \frac{M(0)^{-\gamma}}{\sqrt{\big(-\gamma M(0)^{-\gamma-1}\dot M(0)\big)^2 - 4\gamma^2\mathcal{H}M(0)^{-2\gamma-1}}}.
    \end{equation*}
\end{thm}

\begin{proof}
    Set $N=M^{-\gamma}$, where $\gamma=(k-1)/4$. The virial inequality implies
    \[
        \ddot N \leq \gamma(k+1)\mathcal{H}M^{-\gamma-1}
        =\gamma(k+1)\mathcal{H}N^{1+1/\gamma}.
    \]
    Initially $\dot N(0)<0$. Multiplying by $2\dot N<0$ and integrating gives
    \[
        \dot N(t)^2 \geq \dot N(0)^2+4\gamma^2\mathcal{H}\big(N(t)^{2+1/\gamma}-N(0)^{2+1/\gamma}\big)
        \geq \zeta,
    \]
    where
    \[
        \zeta:=\dot N(0)^2-4\gamma^2\mathcal{H}N(0)^{2+1/\gamma}
        =\gamma^2M(0)^{-2\gamma-2}\big(\dot M(0)^2-4\mathcal{H}M(0)\big)>0.
    \]
    Hence $\dot N(t)\leq-\sqrt{\zeta}$ throughout the forward existence interval. Integrating gives $0\leq N(t)\leq N(0)-\sqrt{\zeta}\,t$, so $N$ reaches zero within time $N(0)/\sqrt{\zeta}$. This yields the stated lifespan bound, with $M(t)\to\infty$ at the maximal forward existence time.
\end{proof}

For the preceding positive-energy theorem, condition $ii)$ and $\dot M^2\leq4ME_k$ imply $E_k>\mathcal H$, hence $E_p=\mathcal H-E_k-E_l<0$.

For the DNKG equation, retaining the linear energy gives another positive-energy criterion.
\begin{thm}[Blow-Up for Positive Energy and Large Displacement]\label{DNKG-large}
    Consider the DNKG with data in $\ell^2\times\ell^2$, $m\ne0$, $\lambda>0$, and $k\geq2$. If
    \[
        \mathcal{H}>0, \qquad \dot M(0)>0, \qquad
        M(0)>\frac{(k+1)\mathcal{H}}{(k-1)m^2},
    \]
    then the solution blows up in finite forward time. In particular,
    \[
        T_+\leq\frac{M(0)}{\gamma\dot M(0)},
        \qquad \gamma=\frac{k-1}{4}.
    \]
\end{thm}

\begin{proof}
    Summation by parts gives
    \[
        E_l=\frac{a}{2}\sumz(q_{n+1}-q_n)^2+m^2M\geq m^2M.
    \]
    The virial identity \eqref{eq:dnkg-virial} holds for all integers $k\geq2$. It gives
    \[
        \ddot M\geq(k-1)m^2M-(k+1)\mathcal H.
    \]
    Set $\delta=(k-1)m^2M(0)-(k+1)\mathcal H>0$. As long as $M\geq M(0)$, the preceding estimate gives $\ddot M\geq\delta$. Since $\dot M(0)>0$, continuity shows that $M$ is increasing throughout its forward existence interval. By the virial identity and $\dot M^2\leq4ME_k$,
    \[
        \ddot N\leq \gamma M^{-\gamma-2}M\big((k+1)\mathcal H-(k-1)m^2M\big)<0,
        \qquad N=M^{-\gamma}.
    \]
    Since $\dot N(0)=-\gamma M(0)^{-\gamma-1}\dot M(0)<0$, concavity gives
    $0<N(t)\leq N(0)+\dot N(0)t$. Thus
    $T_+\leq N(0)/(-\dot N(0))=M(0)/(\gamma\dot M(0))$.
\end{proof}

\section{Examples of Initial Data}

Lastly, let us provide some examples to show that the results above are indeed meaningful. For the FPUT-$\alpha$ system, initial data of the form
\begin{equation*}
    q_n(0) = \begin{cases}
      s & n = 0 \\
      0 & n \neq 0
   \end{cases}, \quad \text{and} \quad p(0) \equiv 0
\end{equation*}
with $|s| < \frac{1}{\sqrt{6}|\alpha|}$ leads to a global solution by Theorem \ref{GLO}. Conversely, let $2 \leq a \in \N$ and $s\in\R$ satisfy $\alpha s > \frac{3}{2}\frac{a}{a-1}$. The ramp data
\begin{equation*}
    q_n(0) = \begin{cases}
      0 & n > a \\
      -\frac{(a-n)}{a} s & 0 \leq n \leq a \\
      0 & n < 0
    \end{cases}, \quad \text{and} \quad p(0) \equiv 0
\end{equation*}
fulfill the criteria of Theorem \ref{E-}, while
\begin{equation*}
    q_n(0) = \begin{cases}
      0 & n > a \\
      -\frac{(a-n)}{a} s & 0 \leq n \leq a \\
      0 & n < 0
   \end{cases}, \quad \text{and} \quad p(0) \equiv c q(0)
\end{equation*}
with the same ramp and sign condition, the data satisfy the conditions of Theorem \ref{E+} for $c >0$ and $c^2>|E_p(0)|/M(0)$, and hence corresponding solutions blow up within finite time. Lastly turning to the DNKG and assuming $a=1$, the conditions of Theorem \ref{DNKG} are fulfilled by the choice
\begin{equation*}
    q_n(0) = \begin{cases} 
      s & n=0 \\
      0 & n \ne 0 
   \end{cases}, \quad s > \bigg(\frac{(2+m^2)(k+1)}{2\lambda}\bigg)^{\frac{1}{k-1}}, \quad \text{and} \quad p(0) \equiv 0.
\end{equation*}

To see that Theorem~\ref{DNKG-large} is not implied by the preceding Theorem~\ref{DNKG-pe}, take $a=1$, $q(0)=s\delta_0$, and $p(0)=cs\delta_0+d\delta_1$, where $\delta_j$ is the standard basis vector at site $j$. Choose
\[
    0<\theta<\frac{(k-1)m^2}{k+1}, \qquad 0<c^2<\theta,
\]
and, for sufficiently large $s>0$, choose $d>0$ so that
\[
    \frac{d^2}{2M(0)}=\theta-c^2-(2+m^2)+\frac{2\lambda}{k+1}s^{k-1}.
\]
Then $\mathcal{H}=\theta M(0)>0$ and $\dot M(0)=2cM(0)>0$. The choice of $\theta$ gives $M(0)>(k+1)\mathcal{H}/((k-1)m^2)$, so Theorem~\ref{DNKG-large} applies. Since $\dot M(0)^2=4c^2M(0)^2<4\mathcal{H}M(0)$, these data fall outside the criterion of Theorem~\ref{DNKG-pe}.

\addcontentsline{toc}{section}{References}

\providecommand{\bysame}{\leavevmode\hbox to3em{\hrulefill}\thinspace}
\providecommand{\MR}{\relax\ifhmode\unskip\space\fi MR }
\providecommand{\MRhref}[2]{%
  \href{http://www.ams.org/mathscinet-getitem?mr=#1}{#2}
}
\providecommand{\href}[2]{#2}

\end{document}